\documentclass[11pt]{article}
\usepackage{amsmath,amsthm,amsfonts}
\usepackage{fullpage}
\usepackage{graphicx}

\newcommand{\C}{\mathbb{C}}
\DeclareMathOperator{\tr}{tr}
\DeclareMathOperator{\normalizer}{\textrm{Norm}}

\DeclareMathAlphabet{\varmathbb}{U}{bbold}{m}{n}
\newcommand{\one}{{\varmathbb 1}}

\begin{document}

\theoremstyle{definition}
\newtheorem{definition}{Definition}
\newtheorem{theorem}{Theorem}
\newtheorem{lemma}{Lemma}
\newtheorem{claim}{Claim}
\newtheorem{proposition}{Proposition}
\newtheorem{exercise}{Exercise}

\newcommand{\Z}{{\mathbb{Z}}}
\newcommand{\Exp}{\mathbb{E}}
\newcommand{\Var}{\mathop\mathrm{Var}}
\newcommand{\Cov}{\mathop\mathrm{Cov}}
\newcommand{\frob}[1]{\left\| #1 \right\|_{\rm frob}}
\newcommand{\opnorm}[1]{\left\| #1 \right\|_{\rm op}}
\newcommand{\norm}[1]{\left\| #1 \right\|}
\newcommand{\abs}[1]{\left| #1 \right|}

\newcommand{\wf}{\widehat{f}}
\newcommand{\wg}{\widehat{g}}

\title{How close can we come to a parity function when there isn't one?}

\author{Cristopher Moore\\Computer Science Department\\University of New Mexico\\and the Santa Fe Institute\\{\tt moore@cs.unm.edu} 
\and Alexander Russell\\Dept. of Computer Science and Engineering\\University of Connecticut\\{\tt acr@cse.uconn.edu}}

\maketitle

\begin{abstract}
Consider a group $G$ such that there is no homomorphism $f:G \to \{ \pm 1 \}$.  In that case, how close can we come to such a homomorphism?  We show that if $f$ has zero expectation, then the probability that $f(xy) = f(x) f(y)$, where $x, y$ are chosen uniformly and independently from $G$, is at most $1/2(1+1/\sqrt{d})$, where $d$ is the dimension of $G$'s smallest nontrivial irreducible representation.  For the alternating group $A_n$, for instance, $d=n-1$.  On the other hand, $A_n$ contains a subgroup isomorphic to $S_{n-2}$, whose parity function we can extend to obtain an $f$ for which this probability is $1/2(1+1/{n \choose 2})$.  Thus the extent to which $f$ can be ``more homomorphic'' than a random function from $A_n$ to $\{ \pm 1 \}$ lies between $O(n^{-1/2})$ and $\Omega(n^{-2})$.
\end{abstract}

The symmetric group $S_n$ has a parity function, i.e., a homomorphism $f: S_n \to \{ \pm 1 \}$, sending even and odd permutations to $+1$ and $-1$ respectively.  The alternating group $A_n$, which consists of the even permutations, has no such homomorphism.  How close can we come to one?  What is the maximum, over all functions $f: A_n \to \{ \pm 1 \}$ with zero expectation, of the probability 
\[
\Pr_{x,y} [f(x)f(y) = f(xy)] 
%= \frac{1}{2} \left( 1 + \Exp_{x,y} \left[ f(x) f(y) f(xy) \right] \right) 
\, ,
\]
where $x$ and $y$ are chosen independently and uniformly from $A_n$?  

We give simple upper and lower bounds on this quantity, for groups in general and for $A_n$ in particular.  Our results are easily extended to functions $f:G \to \C$, but we do not do this here.  Our main result is the following:
\begin{theorem}
\label{thm:upper}
Let $G$ be a group, and let $f:G \to \{ \pm 1 \}$ such that $\Exp f=0$.  Then 
\[
\Pr_{x,y} [f(x)f(y) = f(xy)] \le \frac{1}{2} \left( 1+\frac{1}{\sqrt{d}} \right) \, ,
\]
where $d = \min_{\rho \ne 1} d_\rho$ is the dimension of the smallest nontrivial irreducible representation of $G$.
\end{theorem}
\noindent
Thus if $G$ is \emph{quasirandom} in Gowers' sense~\cite{gowers}---that is, if $\min_{\rho \ne 1} d_\rho$ is large---it is impossible for $f$ to be much more homomorphic than a uniformly random function.  
%, for which $\Pr_{x,y} [f(x)f(y) = f(xy)] = 1/2$.  
For $A_n$ in particular, the dimension of the smallest nontrivial representation is $d=n-1$, so $\Pr_{x,y} [f(x)f(y) = f(xy)] - 1/2 = O(1/\sqrt{n})$.  

If $f$ is a \emph{class function}, i.e., if $f$ is invariant under conjugation so that $f(x^{-1} y x) = f(y)$ for all $x, y \in G$, then we can tighten this bound from $1/\sqrt{d}$ to $1/d$:
\begin{theorem}
\label{thm:upper-class}
Let $G$ be a group, and let $f:G \to \{ \pm 1 \}$ be a class function such that $\Exp f=0$.  Then 
\[
\Pr_{x,y} [f(x)f(y) = f(xy)] \le \frac{1}{2} \left( 1+\frac{1}{d} \right) \, ,
\]
where $d = \min_{\rho \ne 1} d_\rho$ is the dimension of the smallest nontrivial irreducible representation of $G$.
\end{theorem}

As a partial converse to these upper bounds, we have
\begin{theorem}
\label{thm:lower}
Suppose $G$ has a subgroup $H$ with a nontrivial homomorphism $\phi:H \to \{ \pm 1 \}$.  Then there is a function $f:G \to \{ \pm 1 \}$ such that $\Exp f=0$ and 
\[
\Pr_{x,y} [f(x)f(y) = f(xy)] 
\ge \frac{1}{2} \left( 1+ \frac{1}{2} \frac{|H|}{|G|} \left( 1 - \frac{|\normalizer{H}|}{|G|} \right) + \frac{|H|^2}{|G|^2} \right) \, ,
\]
where $\normalizer{H} = \{ c : cHc^{-1} \}$ denotes the normalizer of $H$.
\end{theorem}

If $H$ is normal so that $\normalizer{H}=G$, Theorem~\ref{thm:lower} gives a bias which is quadratically small as a function of the index $|G|/|H|$.  However, in some cases we can do better---for instance, if we can find a set of coset representatives which are involutions:
\begin{theorem}
\label{thm:lower2}
Suppose $G$ has a subgroup $H$ with a nontrivial homomorphism $\phi:H \to \{ \pm 1 \}$.  Suppose further that it has a set of coset representatives $T$ 
%i.e., where $G = \bigcup_{c \in T} cH$ and $|T| = |G|/|H|$, 
such that $c^2=1$ for all $c \in T$.  Then there is a function $f:G \to \{ \pm 1 \}$ such that $\Exp f=0$ and 
\[
\Pr_{x,y} [f(x)f(y) = f(xy)] \ge \frac{1}{2} \left( 1+\frac{|H|}{|G|} \right) \, .
\]
\end{theorem}

For instance, $A_n$ has a subgroup $H$ isomorphic to $S_{n-2}$, consisting of permutations of the last $n-2$ elements, with the first two elements switched if necessary to keep the parity even.  The index of this subgroup is $|H|/|G| = {n \choose 2}$.  Moreover, there is a set of coset representatives $c$ such that $c^2=1$; namely, the permutations which switch the first two elements, setwise, with some other pair.  Thus Theorem~\ref{thm:lower2} applies, and the extent to which $f : A_n \to \{ \pm 1 \}$ can be more homomorphic than a uniformly random function is between $O(n^{-1/2})$ and $\Omega(n^{-2})$.  It would be nice to close this gap.

\begin{proof}[Proof of Theorem~\ref{thm:upper}]
We rely on nonabelian Fourier analysis, for which we refer the reader to~\cite{Serre77}.  In order to establish our notation and choice of normalizations, let $f: G \rightarrow \C$ and let $\rho: G \rightarrow U(d)$ be an irreducible unitary representation of $G$.  We adopt the Fourier transform
$\wf(\rho) = \sum_x f(x) \rho(x)$
in which case we have the Fourier inversion formula
\[
f(x) = \frac{1}{|G|} \sum_\rho d_\rho \tr(\wf(\rho)\;
\rho(x)^\dagger)
\]
and the Plancherel formula
\begin{equation}
\label{eq:inner}
\langle f, g \rangle = \sum_x f(x)^* g(x) = \frac{1}{|G|} \sum_\rho d_\rho
\tr(\wf^\dagger \wg) \, .
\end{equation}
For two functions $f, g: G \rightarrow \C$ we define their convolution
$(f * g)(x) = \sum_y f(y) g(y^{-1}x)$.  With the above normalization,
\[
\widehat{f * g}(\rho) = \wf(\rho) \cdot \wg(\rho) \, .
\]

%\paragraph{Approximate homomorphisms in highly noncommutative groups} 
Now consider a function $f: G \rightarrow \{ \pm 1\}$ 
%$f: G \rightarrow \C$ 
such that 
%$|f(x)| = 1$ and 
$\Exp f = 0$. 
%The condition that $|f(x)| = 1$ implies that $f(x) f(x)^* = 1$. 
We can write the probability that $f$ acts homomorphically on a random pair of elements as an expectation, 
\begin{equation}
\label{eq:as-exp}
\Pr_{x,y} [f(x)f(y) = f(xy)] 
= \frac{1}{2} \left( 1 + \Exp_{x,y} \left[ f(x) f(y) f(xy) \right] \right) \, . 
\end{equation}
We have
\[
%\Exp_{x,y} \left[ f(x) f(y) f(xy)^* \right] 
\Exp_{x,y} \left[ f(x) f(y) f(xy) \right] 
= \Exp_{x,y} \left[ f(x) f(y) g(y^{-1}x^{-1}) \right] = \frac{1}{|G|^2} (f * f * g)(1) \, ,
\]
where 
%$g(z) = f(z^{-1})^*$. 
$g(z) = f(z^{-1})$. 
Observe that
$\wg(\rho) 
%= \sum_x f(x^{-1})^* \rho(x) 
= \sum_x f(x^{-1}) \rho(x) 
%= \sum_x f(x)^* \rho(x)^\dagger 
= \sum_x f(x) \rho(x)^\dagger
= \wf(\rho)^\dagger$ 
and hence, by Fourier inversion, 
\begin{align}
\Exp_{x,y} \left[ f(x) f(y) f(xy) \right] 
&= \frac{1}{|G|^2} (f * f * g)(1) \nonumber \\
%&= \frac{1}{|G|^3} \sum_\rho d_\rho \tr( \wf(\rho) \wf(\rho) \wf(\rho)^\dagger) \nonumber \\
%&= \Exp[f]^3 + \frac{1}{|G|^3} \sum_{\rho \ne 1} d_\rho \tr( \wf(\rho) \wf(\rho) \wf(\rho)^\dagger) \\
&= \frac{1}{|G|^3} \sum_{\rho \ne 1} d_\rho \tr( \wf(\rho) \wf(\rho) \wf(\rho)^\dagger) \, , 
\label{eq:rhorhorho}
\end{align}
where we used the fact that $\wf(1) = |G| \Exp f = 0$.  Everything up to here is essentially identical to the Fourier-analytic treatment of the Blum-Luby-Rubinfeld linearity test~\cite{blum-luby-rubinfeld,bellare-et-al}.

As $N N^\dagger$ is positive semidefinite, 
\begin{equation}
\label{eq:nnn}
\abs{ \tr(N N N^\dagger) } \le \opnorm{N} \tr(N^\dagger N) \le \opnorm{N} \cdot \frob{N}^2 \, , 
\end{equation}
where $\opnorm{N}$ denotes the operator norm 
\[
\opnorm{N} 
%= \max_v \frac{v^\dagger N v}{v^\dagger v} \, , 
= \max_v \frac{\langle v, N v \rangle}{\langle v, v \rangle} \, , 
\]
and $\frob{N}$ denotes the Frobenius norm,
\[
\frob{N} = \sqrt{\tr (N^\dagger N)} = \sqrt{\sum_{ij} \abs{N_{ij}}^2} \, .
\]
Considering also that, from equation~\eqref{eq:inner}, 
\begin{equation}
\label{eq:planch}
\norm{f}^2 = |G| = \langle f, f \rangle = \frac{1}{|G|} \sum_{\rho} d_\rho \frob{\wf(\rho)}^2 \, , 
\end{equation}
we conclude from~\eqref{eq:rhorhorho} and~\eqref{eq:nnn} that
\begin{gather}
\Exp_{x,y} \left[ f(x) f(y) f(xy) \right] 
%- \Exp[f]^3
\le \frac{1}{|G|^3}
\sum_{\rho \ne 1} d_\rho \opnorm{\wf(\rho)} \cdot \frob{\wf(\rho)}^2 
\nonumber \\
\le \max_{\rho \ne 1} \frac{\opnorm{\wf(\rho)}}{|G|^3}
\sum_{\rho \ne 1} d_\rho \frob{\wf(\rho)}^2 
= \max_{\rho \ne 1} \frac{\opnorm{\wf(\rho)}}{|G|} \, .
\label{eq:max-norm}
\end{gather}

Equation~\eqref{eq:planch} also implies that, for any $\rho$, 
%\[
%|G| 
%= \frac{1}{|G|} \sum_\rho d_\rho \frob{\wf(\rho)}^2 \, ,
%%\geq \frac{1}{|G|} \sum_\rho d_\rho \opnorm{\wf(\rho)}^2 \, ,
%%\geq \max_\rho \frac{d_\rho \opnorm{\wf(\rho)}^2}{|G|} \, , 
%\]
\begin{equation}
\label{eq:frob-upper}
\frob{\wf(\rho)} \le \frac{|G|}{\sqrt{d_\rho}} \, . 
\end{equation}
Since $\opnorm{N}$ is $N$'s largest singular value and $\frob{N}^2$ is the sum of their squares, 
\begin{equation}
\label{eq:norm-ineq}
\opnorm{\wf(\rho)} \le \frob{\wf(\rho)} \, .
\end{equation}
Equation~\eqref{eq:frob-upper} then becomes
\begin{equation}
\label{eq:opnorm-upper}
\opnorm{\wf(\rho)} \le \frac{|G|}{\sqrt{d_\rho}} \, . 
\end{equation}
Along with~\eqref{eq:max-norm}, this implies that
\[
\Exp_{x,y} \left[ f(x) f(y) f(xy) \right] 
%- \Exp[f]^3 
\le \max_{\rho \neq 1} \frac{1}{\sqrt{d_\rho}} \, ,
\]
and combining this with~\eqref{eq:as-exp} completes the proof.
\end{proof}

\begin{proof}[Proof of Theorem~\ref{thm:upper-class}]
The proof is the same as that for Theorem~\ref{thm:upper}, except that if $f$ is a class function, then $\wf(\rho)$ is a scalar.  That is, for each $\rho$ there is a $c$ such that $\wf(\rho) = c \one$.  Equation~\eqref{eq:norm-ineq} then becomes
\[
\opnorm{\wf(\rho)} = \abs{c} = \frac{1}{\sqrt{d_\rho}} \frob{\wf(\rho)} \, ,
\]
and~\eqref{eq:opnorm-upper} becomes
\[
\opnorm{\wf(\rho)} \le \frac{|G|}{d_\rho} \, . 
\]
Along with~\eqref{eq:max-norm}, this implies that
\[
\Exp_{x,y} \left[ f(x) f(y) f(xy) \right] 
%- \Exp[f]^3 
\le \max_{\rho \neq 1} \frac{1}{d_\rho} \, ,
\]
and combining this with~\eqref{eq:as-exp} completes the proof as before.
\end{proof}

\begin{proof}[Proof of Theorem~\ref{thm:lower}]
Let $\phi:H \to \{\pm 1\}$ be a homomorphism.  We extend $\phi$ to a function $f:G \to \{\pm 1\}$ in the following way.  We choose a set $T$ of coset representatives such that $G$ is a disjoint union of left cosets, $G = \bigcup_{c \in T} cH$, including the trivial coset $H$ where $c=1$.  Note that $T=|G|/|H|$.  For the trivial coset, we define $f(h) = \phi(h)$ for all $h \in H$.  For each $c \ne 1$, we choose $f(c)$ uniformly from $\{ \pm 1 \}$, and define $f(ch) = f(c) \phi(h)$ for all $h \in H$.  Since $\phi$ is nontrivial, we have $\Exp_H[\phi]=0$ and therefore $\Exp_G[f]=0$.

We will show that, in expectation over $x, y$ and over our choices of $f(c)$, we have 
\begin{equation}
\label{eq:thm-lower-exp}
\Exp[f(x) f(y) f(xy)]
\ge \frac{1}{2} \frac{|H|}{|G|} \left( 1 - \frac{|\normalizer{H}|}{|G|} \right) + \frac{|H|^2}{|G|^2} \, .
\end{equation}
The theorem then follows from~\eqref{eq:as-exp}.

Choose $x, y$ uniformly and independently from $G$.  Write $z=xy$, and consider whether $f(z) = f(x) f(y)$.  There are two cases.  If $y \in H$, then writing $x=ch$ we have
\[
f(z) = f(chy) = f(c) \phi(hy) = f(c) \phi(h) \phi(y) = f(x) f(y) \, . 
\]
The probability of this event is $|H|/|G|$, contributing $|H|/|G|$ to the expectation $\Exp[f(x) f(y) f(xy)]$.

In the other case, $y \in cH$ for some $c \ne 1$.  Then $x$ and $z$ cannot be in the same left coset $c'H$ as each other, since writing $x=c'h$, $y=ck$, and $z=c'\ell$ we would have
\[
c'hck = c'\ell
\]
for some $h,k,\ell \in H$.  This would imply that $hck \in H$ and therefore $c \in H$, a contradiction.  

Now, if $x$ and $z$ are in distinct nontrivial cosets, or if one of $x, z$ is in $H$ but the other is in a nontrivial coset other than $cH$, then $f(x) f(y) f(xy)$ is uniformly random in $\{ \pm 1 \}$.  Thus these events contribute zero to $\Exp[f(x) f(y) f(xy)]$.  This leaves us with two cases: $x \in H$ and $y,z \in cH$, or $x,y \in cH$ and $z \in H$. 

We deal with the case $x \in H$ and $y,z \in cH$ first.  Writing $x=h$, $y=ck$, and $z=c\ell$ gives
\[
hck = c\ell \, ,
\]
or, rearranging, 
\[
c^{-1} h c = \ell k^{-1} \, . 
\]
Then we have 
\[
f(z) = f(c) \phi(\ell) = f(c) \phi(\ell k^{-1}) \phi(k) = f(c) \phi(c^{-1} h c) \phi(k) \, , 
\]
while
\[ 
f(x) f(y) = f(c) \phi(h) \phi(k) \, .
%f(x) f(y) f(xy) = f(c)^2 \phi(h) \phi(k) \phi(\ell) = \phi(h) \phi(k) \phi(\ell) = \phi(h) \phi(\ell k^{-1}) = \phi(h) \phi(c^{-1} h c) \, , 
\]
Thus the question is whether or not
\begin{equation}
\label{eq:conjugate-same}
\phi(c^{-1} h c) = \phi(h) \, .
\end{equation}
The following lemma shows that this is true with probability at least $1/2$ if $h$ is chosen uniformly from $H$ conditioned on $c^{-1} h c \in H$, i.e., uniformly from $H \cap cHc^{-1}$.  Therefore, this event contributes at least zero to $\Exp[f(x) f(y) f(xy)]$.

\begin{lemma}
Let $\phi : H \to \{ \pm 1 \}$ be a homomorphism and let $c \in G$.  Then~\eqref{eq:conjugate-same} holds for at least half the elements of $H \cap cHc^{-1}$.
\end{lemma}

\begin{proof}
We can define a homomorphism $\psi : H \cap cHc^{-1} \to \{ \pm 1 \}$ as 
\[
\psi(h) = \phi(h) \phi(c^{-1} h c) \, .
\]
Clearly~\eqref{eq:conjugate-same} holds if and only if $h \in \ker \psi$, i.e., if $\phi(h)=1$.  But $\ker \psi$ comprises at least half the elements of $H \cap cHc^{-1}$.
\end{proof}

The case $x, y \in cH$ and $z \in H$ is more troublesome.  Writing $x=ch$, $y=ck$, and $z=\ell$, we have
\[
chck=\ell \, .
\]
This event occurs if and only if $chc \in H$.  We then have
\[
f(x) f(y) = f(c)^2 \phi(h) \phi(k) = \phi(h) \phi(k) \, ,
\]
while
\[
f(z) = \phi(\ell) = \phi(chc) \phi(k) \, .
\]
Then, analogous to~\eqref{eq:conjugate-same}, the question is whether
\begin{equation}
\label{eq:chc}
\phi(chc) = \phi(h) \, . 
\end{equation}
Unfortunately, it can be the case that $\phi(chc) = -\phi(c)$ for all $h \in H$ and all $1 \ne c \in T$.  For example, let $G=\{1,c,c^2,c^3\} \cong Z_4$ and $H=\{1,c^2\} \cong \Z_2$, and let $\phi$ be the isomorphism from $H$ to $\{ \pm 1 \}$.  Then $\phi(chc)=-\phi(h)$ for all $h \in H$.

This event, that $chc \in H$ and $\phi(chc) = -\phi(c)$, contributes a negative term to $\Exp[f(x) f(y) f(xy)]$.  We will bound this term by bounding the probability that $chc \in H$ but $c \ne 1$.  
%, where $c \in T$ and $h \in H$ are chosen uniformly at random.  
First consider the following lemma.

\begin{lemma}
\label{lem:chc}
Let $H$ be a subgroup of $G$, let $c \in G$, and suppose that $c \notin \normalizer(H)$.  Then 
\[
\abs{H \cap cHc} \le \frac{|H|}{2} \, . 
\]
\end{lemma}

\begin{proof}
Suppose that $H \cap cHc \ne \emptyset$.  Then there is a pair $h, k \in H$ such that $k = chc$, and
\[
cHc = cHc^{-1} \cdot chc = cHc^{-1} \cdot k \, . 
\]
Since $H = Hk$, we have 
\[
H \cap cHc = ( H \cap cHc^{-1} ) k \, .
\]
In particular, 
\[
\abs{H \cap cHc} = \abs{H \cap cHc^{-1}} \, .
\]
However, if $c \notin \normalizer(H)$ then $H \cap cHc^{-1}$ is a proper subgroup of $H$, in which case its cardinality is at most half that of $H$.
\end{proof}

Now note that if $c' \in H$, then $c'H{c'}^{-1} = cHc^{-1}$.  Therefore, each coset $cH$ is either contained in $\normalizer(H)$ or is disjoint from it.  It follows that the probability that a uniformly random $c \in T$ is in $\normalizer(H)$ is the same as the probability for the entire group, $|\normalizer(H)|/|G|$.   Since $|T| = |G|/|H|$, the number of such $c$ is
\[
\abs{T \cap \normalizer(H)} = \frac{|\normalizer(H)|}{|H|} \, . 
\]
Thus we have $|\normalizer(H)/|H|-1$ coset representatives $c \in \normalizer(H)$ other than $c=1$.  If we condition on the event that $x, y \in cH$, each of these $c$ could conceivably contribute $-1$ to $\Exp[f(x) f(y) f(xy)]$, while Lemma~\ref{lem:chc} implies that the other $|G|/|H| - |\normalizer(H)/|H|$ coset representatives contribute at least $-1/2$.  
% to $\Exp[f(x) f(y) f(xy)]$.  
The total contribution of the case $x,y \in cH$, $z \in H$ to $\Exp[f(x) f(y) f(xy)]$ is then at least
\begin{align*}
- &\frac{|H|^2}{|G|^2} \left( \frac{|\normalizer(H)|}{|H|} - 1 + \frac{1}{2} \left( \frac{|G|}{|H|} - \frac{|\normalizer(H)|}{|H|} \right) \right) \\
= &\frac{1}{2} \frac{|H|}{|G|} \left( -1 - \frac{|\normalizer(H)|}{|G|} \right) 
+ \frac{|H|^2}{|G|^2} \, .
\end{align*}
Adding the contribution $|H|/|G|$ from the case $y \in H$ gives~\eqref{eq:thm-lower-exp} and completes the proof.
\end{proof}

\begin{proof}[Proof of Theorem~\ref{thm:lower2}]
If $c^2=1$, then $cHc=cHc^{-1}$.  This changes the troublesome case to the easy one, where~\eqref{eq:conjugate-same} holds with probability at least $1/2$ for all $h \in H \cap cHc^{-1}$, and so the event $x, y \in cH$, $z \in H$ contributes at least zero to $\Exp[f(x) f(y) f(xy)]$.  We then have $\Exp[f(x) f(y) f(xy)] \ge |H|/|G|$ from the case $y \in H$, and the theorem follows from~\eqref{eq:as-exp}.
\end{proof}

Note that the premise of Theorem~\ref{thm:lower2} can be weakened considerably: namely, that for all $c$ such that $H \cap cHc \ne \emptyset$, we have $c^2 = k$ for some $k \in H$ with $\phi(k)=1$.

\paragraph{Acknowledgments}  We benefited from discussions with Avi
Wigderson, who obtained Theorem~\ref{thm:upper} and generalizations
independently, and from his lectures and those of Ben Green at the
Bellairs Research Institute of McGill University.  This work was
supported by the NSF under grants CCF-0829931, 0835735, and 0829917 and by the DTO under contract W911NF-04-R-0009.

\end{document}